\documentclass[a4paper, draft]{article}

\usepackage[leqno]{amsmath}
\usepackage{amssymb} 
\usepackage{amsthm}  
\usepackage{hhline}
\usepackage{array}
\usepackage{afterpage}
\usepackage{doublespace}

\newtheorem{tm}{Theorem}[section]

\newtheorem{pr}[tm]{Proposition}

\newcommand{\fej}[1]{\section{\!\!\!\!\!\!.\ #1}}
\newcommand{\alfej}[1]{\subsection{\!\!\!#1}}
\newcommand{\nemfej}[1]{\renewcommand\thesection{}\section{\!\!\!\!\!\!#1}\renewcommand\thesection{\arabic{section}}}

\newcommand*{\hop}{\bigskip\\}
\renewcommand*{\int}{\intop\limits}
\newcommand*{\Zb}{\mathbb Z}
\newcommand*{\Rb}{\mathbb R}

\newcommand*{\Nb}{\mathbb N}
\newcommand*{\om}{\omega}

\newcommand*{\un}[1]{\underline{#1}}

\newcommand*{\be}{\beta}

\newcommand*{\e}[1]{\text{\rm e}^{#1}}

\newcommand*{\ve}{\varepsilon}

\newcommand*{\di}{\,\text{\rm d}}

\newcommand*{\Ev}{{\bf E}}

\newcommand*{\vp}{\varphi}
\newcommand*{\tv}{\mbox{\boldmath$\tau$}}
\newcommand*{\pt}{\partial}
\newcommand*{\te}{\theta}

\DeclareMathOperator{\ctanh}{ctanh}

\begin{document}
\begin{spacing}{1.3}
\title{Microscopic shape of shocks in a domain growth model}
\author{M\'arton Bal\'azs\\
\medskip
Institute of Mathematics, TU Budapest}
\date{}
\maketitle

\begin{abstract}
Considering the hydrodynamical limit of some interacting particle systems leads to hyperbolic differential equation for the conserved quantities, e.g.\! the inviscid Burgers equation for the simple exclusion process. The physical solutions of these partial differential equations develop discontinuities, called shocks. The microscopic structure of these shocks is of much interest and far from being well understood. We introduce a domain growth model in which we find a stationary (in time) product measure for the model, as seen from a defect tracer or second class particle, traveling with the shock. We also show that under some natural assumptions valid for a wider class of domain growth models, no other model has stationary product measure as seen from the moving defect tracer.
\end{abstract}
\bigskip
{\bf Key-words:} second class particle; shock solution.

\nemfej{Introduction}

The hydrodynamical limit of the nearest neighbor asymmetric simple exclusion model leads to the inviscid Burgers equation
\[
\frac{\pt u}{\pt t}+\frac12\,\frac{\pt u^2}{\pt x}=0
\]
which is a special case of the one-component hyperbolic conservation law
\begin{equation}
\frac{\pt u}{\pt t}+\frac{\pt J(u)}{\pt x}=0\label{eq:hydalap}
\end{equation}
where $u\mapsto J(u)$ is a smooth, typically convex function. (By changing $x$ to $-x$, concave $J$-s can be transformed to convex ones.) This equation has a shock (weak) solution starting with initial data 
\[
u(0,\,x)=\left\{\begin{array}{lcl}
\!\!\!u_{\text{left}}&,\ \ &x<0\\
\!\!\!u_{\text{right}}&,\ \ &x\ge0
\end{array}\right. 
\]
with $u_{\text{left}}>u_{\text{right}}$. The stable weak solution is of the form
\[
u(t,\,x)=\left\{\begin{array}{lcl}
\!\!\!u_{\text{left}}&,\ \ &x<st\\
\!\!\!u_{\text{right}}&,\ \ &x\ge st
\end{array}\right. 
\]
where the speed $s$ of the traveling shock is determined by the Rankine-Hugoniot formula
\begin{equation}
s=\frac{J(u_{\text{right}})-J(u_{\text{left}})}{u_{\text{right}}-u_{\text{left}}}\ \ ,\label{eq:ran}
\end{equation}
see e.g.\ \cite{smo}. This is what we see on a macroscopic scale. The microscopic structure (i.e.\ on the level of particles) of the shock is of great interest. It has been considered in the context of the asymmetric simple exclusion process, and rather complicated microscopic structures have been found \cite{mkps} \cite{dls} \cite{shock} \cite{imes} \cite{sps}. In the more general context of attractive particle systems the microscopic structure of the shock was investigated by \cite{rez}. 

In the present note we consider a class of one-dimensional domain growth models, pa\-ra\-met\-ri\-sed by a jump rate function, $r\,:\,\Zb\to\Rb$. In a special case of the rate function we show that the shock, as seen from a defect tracer (second class particle) has stationary (in time) distribution of product structure which we identify. We also show that this is a peculiarity of the case considered, no other model in the wide class of these models has this property. The structure of the paper is the following: 

\noindent In section \ref{sc:model} we define the class of models considered and determine the stationary distributions for them. 

\noindent We describe the hydrodynamic limit of these models and calculate the speed of the shocks using Rankine-Hugoniot formula \eqref{eq:ran} in section \ref{sc:hyd}.

\noindent In section \ref{sc:coupl} we introduce the defect tracer in our models. Via Rankine-Hugoniot formula, we also give an indication on the fact that, in general, shock solutions are closely related to measures stationary as seen from the defect tracer.

\noindent The last section contains our main result on the product structure of such a stationary distribution as seen from the defect tracer. This gives an explicit description of the microscopic shape of some types of shock solutions.

\fej{The bricklayers' model}\label{sc:model}

\alfej{Infinitesimal generator}

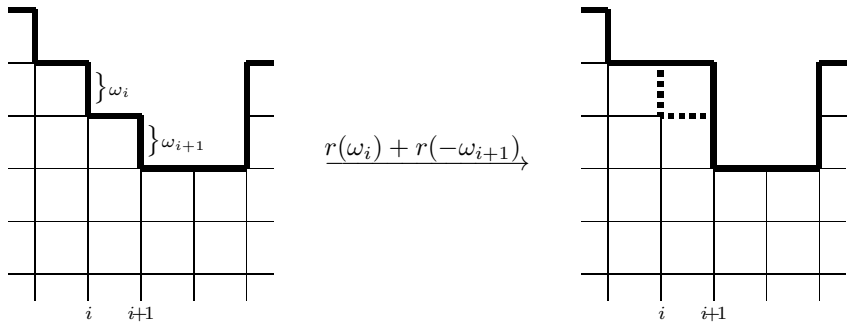
\begin{figure}[p]
\begin{center}
\begin{picture}(100, 160)(0, -10)
\linethickness{0.2pt}
\put(0, 10){\line(1, 0){100}}
\put(0, 30){\line(1, 0){100}}
\put(0, 50){\line(1, 0){100}}
\put(0, 70){\line(1, 0){30}}
\put(90, 70){\line(1, 0){10}}
\put(0, 90){\line(1, 0){10}}

\put(10, 0){\line(0, 1){90}}
\put(30, 0){\line(0, 1){70}}
\put(50, 0){\line(0, 1){50}}
\put(70, 0){\line(0, 1){50}}
\put(90, 0){\line(0, 1){50}}

\put(29, -7){$\scriptstyle{i}$}
\put(45, -7){$\scriptstyle{i\!+\!1}$}

\linethickness{2pt}

\put(0, 110){\line(1, 0){10}}
\put(10, 90){\line(1, 0){20}}
\put(30, 70){\line(1, 0){20}}
\put(50, 50){\line(1, 0){40}}
\put(90, 90){\line(1, 0){10}}

\put(10, 90){\line(0, 1){20}}
\put(30, 70){\line(0, 1){20}}
\put(50, 50){\line(0, 1){20}}
\put(90, 50){\line(0, 1){40}}

\put(32, 78){$\scriptstyle{\bigl\}\om_i}$}
\put(52, 58){$\scriptstyle{\bigl\}\om_{i+1}}$}

\end{picture}
\begin{picture}(110, 160)
\put(30, 60){\makebox(50, 12){$\underrightarrow{r(\om_i)+r(-\om_{i+1})\ }$}}
\end{picture}
\begin{picture}(100, 160)(0, -10)

\linethickness{0.2pt}

\put(0, 10){\line(1, 0){100}}
\put(0, 30){\line(1, 0){100}}
\put(0, 50){\line(1, 0){100}}
\put(0, 70){\line(1, 0){30}}
\put(90, 70){\line(1, 0){10}}
\put(0, 90){\line(1, 0){10}}

\put(10, 0){\line(0, 1){90}}
\put(30, 0){\line(0, 1){70}}
\put(50, 0){\line(0, 1){50}}
\put(70, 0){\line(0, 1){50}}
\put(90, 0){\line(0, 1){50}}

\put(29, -7){$\scriptstyle{i}$}
\put(45, -7){$\scriptstyle{i\!+\!1}$}

\linethickness{2pt}

\put(0, 110){\line(1, 0){10}}
\put(10, 90){\line(1, 0){40}}
\put(50, 50){\line(1, 0){40}}
\put(90, 90){\line(1, 0){10}}

\put(10, 90){\line(0, 1){20}}
\put(50, 50){\line(0, 1){40}}
\put(90, 50){\line(0, 1){40}}

\put(30, 70){\dashbox{2}(20, 20){}}

\end{picture}
\end{center}
\caption{A possible move}\label{fig:elso}
\end{figure}

We consider the phase space 
\begin{equation*}
\Omega=\left\{\un\om=(\om_i)_{i\in\Zb}\ :\ \om_i\in\Zb\right\}=\Zb^{\Zb}\ \ .
\end{equation*}
For each pair of neighboring sites $i$ and $i+1$ of $\Zb$, we can imagine a column built of bricks, above the edge $(i,\,i+1)$. The height of this column is denoted by $h_i$. If $\un{\om}(t)\in\Omega$ for a fixed time $t\in\Rb$ then $\om_i(t)=h_{i-1}(t)-h_i(t)\,\in\Zb$ is the negative discrete gradient of the height of the ``wall''. The growth of a column is described by Poisson processes. A brick can be added to a column:
\[
\left(\om_i,\,\om_{i+1}\right)\longrightarrow\left(\om_i-1,\,\om_{i+1}+1\right)\ \ \text{with rate}\ \ r(\om_i)+r(-\om_{i+1})\ \ .
\]
See fig.\,\ref{fig:elso} for some possible instantaneous changes. The process can be represented by bricklayers standing at each site $i$, laying a brick on the column on their right with rate $r(\om_i)$ and laying a brick to their left with rate $r(-\om_i)$. This interpretation gives reason to call these model bricklayers' model. %
\afterpage{\clearpage}%
For small $\ve$ the conditional expectation of the growth of the column between $i$ and $i+1$ in the time interval $[t,\,t+\varepsilon]$ is $\left\{r(\om_i(t))+r(-\om_{i+1}(t))\right\}\cdot\varepsilon+\mathfrak{o}(\ve)$. Note that the process has a left-right mirror symmetry, i.e.\ the rate of a column's growth is the same as if looking at the reflected configuration. We want the dynamics to smoothen our interface, that is why we assume monotonicity of the rate function $r$, which means that a column grows more rapidly if it has a higher neighbor on the right or on the left. In later sections we shall impose another restrictive condition on $r$, see \eqref{eq:ratafelt}.

At time $t$, the interface mentioned before is described by $\un{\om}(t)$. Let $\varphi\,:\,\Omega\to\Rb$ be a bounded cylinder function i.e.\ $\varphi$ depends on a finite number of values of $\om_i$. The growth of this interface is a Markov process, with the formal infinitesimal generator $L$:
\[
(L\varphi)(\un\om)=\sum_{i\in\Zb}\Bigl\{\left[r(\om_i)+r(-\om_{i+1})\right]\cdot\left[\varphi(\dots,\,\om_i-1,\,\om_{i+1}+1,\,\dots)-\varphi(\un\om)\right]\Bigr\}\ \ .
\]
Note that for each index $i,\ \om_i$ can also be negative hence direct particle interpretation fails, see the remark after formula \eqref{eq:kons}. 

When constructing the process rigorously, problems may arise due to the unbounded growth rates. The system being one-component and attractive, we assume that existence of dynamics on a set of tempered configurations $\widetilde\Omega$ (i.e.\ configurations obeying some restrictive growth conditions) can be established by applying methods initiated by Liggett and Andjel \cite{and} \cite{exi}. Technically we assume that $\widetilde\Omega$ is of full measure w.r.t.\ the canonical Gibbs measures defined in \ref{sc:gibbs}. We do not deal with this question in the present paper.
\hop
\noindent {\bf The exponential bricklayers' model}

\bigskip

\noindent A special case of the models is the exponential bricklayers' model (EBL), where for $z\in\Zb$
\begin{equation}
r(z)=\e{-\frac{\be}{2}}\,\e{\be z}\label{eq:eblr}
\end{equation}
with a positive real parameter $\be$.

\alfej{Translation invariant stationary product measures}\label{sc:gibbs}

In this subsection we show a natural way to construct a stationary translation invariant product measure for our models. By chapter one of \cite{ips}, a measure $\mu$ is stationary, iff for any bounded cylinder function $\vp$,
\[
\Ev(L\,\vp)(\un\om)=0
\]
is satisfied for a process distributed according to $\mu$. We assume $\mu$ to be a product measure with marginals 
\[
\mu(z)=\mu\{\un\om\,:\,\om_i=z\}
\]
for $z\in\Zb$. By changing variables and using product structure of $\mu$,
\begin{multline*}
\Ev(L\,\vp)(\un\om)=\\
=\Ev\sum_{i\in\Zb}\Bigl\{\left[r(\om_i)+r(-\om_{i+1})\right]\cdot\left[\varphi(\dots,\,\om_i-1,\,\om_{i+1}+1,\,\dots)-\varphi(\un\om)\right]\Bigr\}=\\
=\Ev\sum_{i\in\Zb}\biggl[r(\om_i+1)\cdot\frac{\mu(\om_i+1)}{\mu(\om_i)}\cdot\frac{\mu(\om_{i+1}-1)}{\mu(\om_{i+1})}+r(-\om_{i+1}+1)\cdot\frac{\mu(\om_i+1)}{\mu(\om_i)}\cdot\frac{\mu(\om_{i+1}-1)}{\mu(\om_{i+1})}\\
-r(\om_i)-r(-\om_{i+1})\biggr]\cdot\varphi(\un\om)\ \ .
\end{multline*}
This expression becomes zero if we make the sum telescopic on the cylinder set supporting $\vp$. Hence stationarity of $\mu$ is assured by assuming
\[
r(z)\cdot\frac{\mu(z)}{\mu(z-1)}\ \ \text{and}\ \ r(-z)\cdot\frac{\mu(z)}{\mu(z+1)}
\]
to be constants. As a consequence, we obtain the condition
\begin{equation}
r(z)\cdot r(-z+1)=\text{constant\ \ .}\label{eq:kons}
\end{equation}
There are two essentially different choices. 
\[
r(z)\cdot r(-z+1)=0
\]
defines models of zero range types, we do not consider this possibility here. The other choice is choosing the right-hand side of \eqref{eq:kons} to be a positive constant. In this case, by rescaling time, we can turn this constant to be one without loss of generality:
\begin{equation}
r(z)\cdot r(-z+1)=1\ \ .\label{eq:ratafelt}
\end{equation}
Rates \eqref{eq:eblr} of the EBL model satisfy this condition. 

For $n\in\Nb$, we define 
\[
r(n)!:\,=\prod_{y=1}^nr(y)
\]
with the convention that the empty product has value one. Let
\[
\bar\te:\,=\log\left(\liminf_{n\to\infty}\left(r(n)!\right)^{1/n}\right)=\lim_{n\to\infty}\log(r(n))\ \ ,
\]
which is strictly positive by \eqref{eq:ratafelt} and by monotonicity of $r$, and can even be infinite. With a generic real parameter $\te\in\left(-\bar\te,\,\bar\te\right)$, we define 
\[
Z(\te):\,=\sum_{z=-\infty}^\infty\frac{\e{\te z}}{r(|z|)!}
\]
and the product measure $\mu^{(\te)}$ with marginals
\begin{equation}
\mu^{(\te)}(z):\,=\frac{1}{Z(\te)}\cdot\frac{\e{\te z}}{r(|z|)!}\ \ ,\label{eq:om}
\end{equation}
which has the property
\begin{equation}
\begin{split}
r(z)\cdot\frac{\mu^{(\te)}(z)}{\mu^{(\te)}(z-1)}&=\e{\te}\\
r(-z)\cdot\frac{\mu^{(\te)}(z)}{\mu^{(\te)}(z+1)}&=\e{-\te}\ \ ,
\end{split}\label{eq:miap}
\end{equation}
thus it is stationary. We call these measures canonical Gibbs-measures.

For the special case of the EBL model, for $\te\in(-\infty,\,\infty)$, we obtain the discrete normal distribution
\begin{equation}
\mu^{(\te)}(z)=\frac{\e{-\frac{\be}{2}\left(z-\frac{\te}{\be}\right)^2}}{\e{-\frac{\te^2}{2\be}}Z(\te)}=\frac{\e{-\frac{\be}{2}(z-m)^2}}{\widetilde{Z}(\be,\,m)}\label{eq:origmu}
\end{equation}
with the notation $m:=\te/\be$.

\fej{Hydrodynamical limit}\label{sc:hyd}

Being attractive due to monotonicity of $r$, we can take the hydrodynamical limit of a bricklayers' model by 
\begin{equation}
u(t,\,x):\,=\Ev\,\om_{x/\ve}(t/\ve)\ \ .\label{eq:miau}
\end{equation}
Then via formal computations we obtain differential equation \eqref{eq:hydalap}
\[
\frac{\pt u}{\pt t}+\frac{\pt J(u)}{\pt x}=0
\]
as $\ve\to0$, where $J(u)$ is defined as follows. The function $u(\te)=\Ev^{(\te)}(\om)$ of $\te$ is strictly increasing since the derivative
\[
\frac{\di u(\te)}{\di\te}=\left(\Ev^{(\te)}(\om^2)-\left(\Ev^{(\te)}(\om)\right)^2\right)
\]
is positive ($-\bar\te<\te<\bar\te$). Let $\te(u)$ be the inverse function. The quantity $\Ev^{(\te)}\left(r(\om)+r(-\om)\right)$ depends on $\te$, and $J$ is defined by
\begin{equation}
J(u):=\Ev^{(\te(u))}\left(r(\om)+r(-\om)\right)=2\cosh(\te(u))\ \ .\label{eq:miaj}
\end{equation}
\begin{pr}
There exist $\te_1<0<\te_2$ numbers such that $J(u)$ defined above is convex on the interval $\left(u(\te_1),\,u(\te_2)\right)$.
\end{pr}
\begin{proof}
With the notations
\[
u'(\te):\,=\frac{\di u(\te)}{\di\te}\ \ \text{and}\ \ u''(\te):\,=\frac{\di^2u(\te)}{\di\te^2}
\]
and by computing derivatives of inverse functions, we obtain from \eqref{eq:miaj}
\[
\frac{\di^2 J}{\di u^2}\circ(u(\te))=\frac{\cosh(\te)}{2}\,\frac{1}{(u'(\te))^2}-\frac{\sinh(\te)}{2}\,\frac{u''(\te)}{(u'(\te))^3}\ \ .
\]
The positivity of the left-hand side is assured in case $\te=0$, and is equivalent to the condition
{\renewcommand{\arraystretch}{2.2}
\begin{equation}
\begin{array}{ll}
\displaystyle\frac{u''(\te)}{u'(\te)}<\ctanh(\te)\ ,\ &\text{if}\ \te>0\ \ ,\\
\displaystyle\frac{u''(\te)}{u'(\te)}>\ctanh(\te)\ ,\ &\text{if}\ \te<0
\end{array}\label{eq:convfelt}
\end{equation}
}by positivity of $u'(\te)$. The function $\te\mapsto u(\te)$ is analytic in $(-\bar\te,\,\bar\te),\ u'(\te)$ is strictly positive, hence the left-hand side of \eqref{eq:convfelt} is bounded on the interval $(-\te^*,\,\te^*)$ for each $0<\te^*<\bar\te$. Due to the unbounded behavior of $\ctanh(\te)$ on any interval containing zero, there exist $\te_1<0<\te_2$ for which \eqref{eq:convfelt} and hence convexity of $J(u)$ is satisfied.
\end{proof}

Using definitions \eqref{eq:miaj} and \eqref{eq:miau} in Rankine-Hugoniot formula \eqref{eq:ran}, the speed of the shock can now be written as
\begin{equation}
s=\frac{\Ev^{(\te(u_{\text{right}}))}\left(r(\om)+r(-\om)\right)-\Ev^{(\te(u_{\text{left}}))}\left(r(\om)+r(-\om)\right)}{\Ev^{(\te(u_{\text{right}}))}(\om)-\Ev^{(\te(u_{\text{left}}))}(\om)}\ \ .\label{eq:mitran}
\end{equation}

\fej{The defect tracer}\label{sc:coupl}

\alfej{Coupling the models}

Let $\un\om^+(0)$ and $\un\om^-(0)$ be two elements of $\widetilde\Omega$. At time $t=0$ we start with a configuration where these two realizations differ at only one site:
\[
\om^+_i(0)=\om^-_i(0)\ \ \text{if\ \ } i\ne0\ \ \ \text{and}\ \ \ \om^+_0(0)=\om^-_0(0)+1\ \ .
\]
One possible representation can be imagined by two walls. At time 0, the walls are the same on the right side of position 0, and every column of the wall $^+$ is higher by one brick than column of wall $^-$ on the left side of zero. We want the two processes to grow together in such a way, that the difference between them remains ``one step'' at any time $t>0$:
\begin{multline*}
(\forall\,t>0)\ (\exists_1\,Q(t)\in\Zb) : \om^+_i(t)=\om^-_i(t)\ \ \text{if\ \ } i\ne Q(t)\ \ \ \text{and}\\
\om^+_{Q(t)}(t)=\om^-_{Q(t)}(t)+1\ \ .
\end{multline*}
We shall call this difference between the two models \emph{defect tracer}, and $Q(t)$ is its position at time $t$. We show the coupling which preserves the only one defect tracer while both $\un\om^-$ and $\un\om^+$ evolves as usual. This coupling for the simple exclusion model is described (with particle notations) in \cite{couse} and \cite{ips}. Let our defect tracer be posed at point $Q$ (i.e.\ $\om_Q^+=\om_Q^-+1;\ \om_i^+=\om_i^-$\ if $i\ne Q$), and let $h^+_i\uparrow$ (or $h^-_i\uparrow$) mean that the column of $\un\om^+$ (or the column of $\un\om^-$, respectively) between the points $i$ and $i+1$ has grown by one brick. Then the growing rule for the columns $h^\pm_{Q-1}$ and $h^\pm_Q$ is shown in table \ref{tab:cou}.%
\renewcommand{\arraystretch}{1.4} 
\begin{table}[p]
\[
\begin{array}{|c||c|c|c|c||c|}
\hline
\text{with rate}&h_{Q-1}^-\uparrow&h_{Q-1}^+\uparrow&h_Q^-\uparrow&h_Q^+\uparrow&Q\ \text{has\dots}\\
\hhline{|=#=|=|=|=#=|}
r(-\om_Q^-)-r(-\om_Q^+)&\bullet&&&&\text{decreased}\\
\hline
r(\om_{Q-1}^-)+r(-\om_Q^+)&\bullet&\bullet&&&-\\
\hline
r(\om_Q^+)-r(\om_Q^-)&&&&\bullet&\text{increased}\\
\hline
r(\om_Q^-)+r(-\om_{Q+1}^-)&&&\bullet&\bullet&-\\
\hhline{|=#=|=|=|=#=|}
\end{array} 
\]
\caption{The coupling rules\label{tab:cou}}
\end{table}\afterpage{\clearpage}%
Every line of that table represents a possible step with rate written on the first column. These rates are non-negative due to monotonicity of $r$. For each column of this table, summing the rates corresponding to the possible steps assures us that columns of each $\un\om^+$ and $\un\om^-$ evolve as usual in the neighborhood of $Q$. For $i\ne Q-1$ or $Q,\ \ h_i^+$ and $h_i^-$ increases at the same time with the original rate $r(\om_i^-)+r(-\om_{i+1}^-)=r(\om_i^+)+r(-\om_{i+1}^+)$. 

How does an observer following the defect tracer see the surface? We introduce the drifted form $\tv_k\,\un\om$ of an $\un\om\in\Omega$ as follows. Let $k\in\Zb$, then $\tv_k\,\un\om\in\Omega$ and
\[
(\tv_k\,\un\om)_i:\,=\om_{i-k}\ \ .
\]

From now on, we denote by $\un\om(t)$ the $\un\om^-(t)$ process as seen from the position $Q(t)$ of the defect tracer, i.e.\ $\om_i:\,=\om_{Q+i}^-$. According to the coupling rules, we can write the infinitesimal generator for $\un\om$:
\begin{multline}
(L^{\text{(d.t.)}}\varphi)(\un\om)=\\
=\sum_{i\ne-1,\,0}\Bigl\{\left[r(\om_i)+r(-\om_{i+1})\right]\cdot\left[\varphi(\dots,\,\om_i-1,\,\om_{i+1}+1,\,\dots)-\varphi(\un\om)\right]\Bigr\}+\\
+\left[r(\om_{-1})+r(-\om_0-1)\right]\cdot\left[\vp(\dots,\,\om_{-1}-1,\,\om_0+1,\,\dots)-\vp(\un\om)\right]+\\
+\left[r(\om_0)+r(-\om_1)\right]\cdot\left[\varphi(\dots,\,\om_0-1,\,\om_1+1,\,\dots)-\vp(\un\om)\right]+\\
+\left[r(-\om_0)-r(-\om_0-1)\right]\cdot\left[\vp(\tv_1(\dots,\,\om_{-1}-1,\,\om_0+1,\,\dots))-\vp(\un\om)\right]+\\
+\left[r(\om_0+1)-r(\om_0)\right]\cdot\left[\vp(\tv_{-1}\,\un\om)-\vp(\un\om)\right]\ \ .\label{eq:utgen}
\end{multline}

\alfej{The speed of the defect tracer}

The main problem of this note is to find a stationary measure for the process as seen from the defect tracer, i.e.\ to find a measure $\mu^{\text{(d.t.)}}$, for which
\[
\Ev(L^{\text{(d.t.)}}\vp)(\un\om)=0
\]
is satisfied. Before giving a partial answer to this question, we give an early indication on the correspondence to shocks of such a measure $\mu^{\text{(d.t.)}}$. 

Let $a<-1$ (and $b>1$) be sites far on the left side (and far on the right side, respectively) of the defect tracer. We choose the function 
\[
\vp(\un\om):=\sum_{k=a}^b\om_k
\]
in \eqref{eq:utgen} to obtain
\begin{multline*}
(L^{\text{(d.t.)}}\vp)(\un\om)=\left[r(\om_{a-1})+r(-\om_a)\right]-\left[r(\om_b)+r(-\om_{b+1})\right]+\\
+\left[r(-\om_0)-r(-\om_0-1)\right]\cdot(\om_{a-1}-\om_b)+\left[r(\om_0+1)-r(\om_0)\right]\cdot(\om_{b+1}-\om_a)\ \ .
\end{multline*}
Let us assume that a measure $\mu^{\text(d.t.)}$ is stationary for $L^{\text(d.t.)}$. Let us also assume that as $l\to\pm\infty$, the random variable $\om_l$ becomes asymptotically independent of $\om_{-1},\,\om_0,\,\om_1$, and the distribution of $\om_l$ converges weakly. Then we have
\begin{multline}
0=\Ev(L^{\text{(d.t.)}}\vp)(\un\om)=
\Ev\left[r(\om_a)+r(-\om_a)\right]-\Ev\left[r(\om_b)+r(-\om_b)\right]+\\
+\Ev\left[r(-\om_0)-r(-\om_0-1)\right]\cdot\Ev(\om_a-\om_b)+\Ev\left[r(\om_0+1)-r(\om_0)\right]\cdot\Ev(\om_b-\om_a)+\\
+H(a,\,b)\ \ ,\label{eq:diszran}
\end{multline}
where the error function $H(a,\,b)$ tends to zero if $a\to-\infty$ and $b\to\infty$. (For the product measure $\mu^{\text{(d.t.)}}$ we find in the next section, $H(a,\,b)=0$ for any $a<-1,\ b>1$.) According to the rules of the coupling, and assuming also ergodicity of the process as seen from the view of the defect tracer, we have the law of large numbers
\[
v:=\lim_{t\to\infty}\frac{Q(t)}{t}=\Ev\left\{\left[r(\om_0+1)-r(\om_0)\right]-\left[r(-\om_0)-r(-\om_0-1)\right]\right\}\ \ \text{a.s.}
\]
for the speed of the defect tracer. Hence we conclude from \eqref{eq:diszran} that
\[
v=\lim_{a\to-\infty,\,b\to\infty}\frac{\Ev\left[r(\om_b)+r(-\om_b)\right]-\Ev\left[r(\om_a)+r(-\om_a)\right]}{\Ev(\om_b)-\Ev(\om_a)}
\]
in case $\Ev(\om_a)\ne\Ev(\om_b)$ and their limits are not equal i.e.\ the slope of the surface is different far on the two sides. This formula is the same as \eqref{eq:mitran}, which we obtained for the speed of the shock using the Rankine-Hugoniot formula. This shows that a measure $\mu^{\text{(d.t.)}}$ with different asymptotics on the left and on the right can be identified as the microscopic structure of a shock solution of \eqref{eq:hydalap}.

\fej{Stationary measures as seen from the defect tracer}

In this section we find a stationary product measure satisfying \eqref{eq:miap} for the defect tracer of the EBL model. We also show that this kind of measure only exists for the EBL model. 

Intuitively one expects that far from the defect tracer a stationary measure behaves like the canonical Gibbs-measure $\mu^{(\te)}$, since the defect tracer is only a local ``error'' for the evolution of the process. The canonical measure has one parameter $\te$, but it is not necessary that in this case parameter $\te_{\text{left}}$ far on the left side would be equal to the parameter $\te_{\text{right}}$ far on the right side. Let 
\[
\un\te:\,=\{\te_i\,:\,i\in\Zb\}
\]
be a sequence of parameters. Then it seems to be reasonable to assume that the product measure $\mu^{(\un\te)}$ with marginals
\[
\mu_i(z)=\mu^{(\un\te)}\left\{\un\om\,:\,\om_i=z\right\}:\,=\mu^{(\te_i)}(z)=\frac{1}{Z(\te_i)}\cdot\frac{\e{\te_i z}}{r(|z|)!}
\]
is stationary for $L^{\text{d.t.}}$ (\ref{eq:utgen}) ($i\in\Zb$). This measure only differs from the canonical $\mu^{(\te)}$ \eqref{eq:om} in that the parameter of its one-dimensional marginals depends on the position. The question is whether there are any choices of $\un\te$ for $\mu^{(\un\te)}$ to be stationary.
\begin{tm}
For a bricklayers' model, if $r$ is not the constant function, then the measure $\mu^{(\un\te)}$ described above is stationary for $L^{\text{(d.t.)}}$ if and only if $r$ is the rate of an EBL model with any parameter $\be>0$, and for the $\un\te$ parameters of $\mu^{(\un\te)}$
{\renewcommand{\arraystretch}{1}
\[
\te_i=\left\{\begin{array}{ll}\te_{\text{left}}&\text{if}\ \ i\le-1\ ,\\\te_{\text{right}}:\,=\te_{\text{left}}-\be\ \ &\text{if}\ \ i\ge0\end{array}\right.
\]
}is satisfied with an arbitrary real number $\te_{\text{left}}$.
\end{tm}
\begin{proof}
Stationarity means 
\[
\Ev^{(\un\te)}(L^{\text{d.t.}}\vp)(\un\om)=0\ \ ,
\]
after some changes of variables, by straightforward computations we obtain from (\ref{eq:utgen})
\begin{equation}
0=\Ev^{(\un\te)}\Biggl\{\Bigl\{A+B+C+D\Bigr\}\vp(\un\om)\Biggr\}\ \ ,\label{eq:foell}
\end{equation}
where
\begin{eqnarray*}
A&=&\sum_{i\ne-1}\Bigl[\left[r(\om_i+1)+r(-\om_{i+1}+1)\right]\cdot\frac{\mu_i(\om_i+1)\,\mu_{i+1}(\om_{i+1}-1)}{\mu_i(\om_i)\,\mu_{i+1}(\om_{i+1})}-\\
&-&\left[r(\om_i)+r(-\om_{i+1})\right]\Bigr]\ \ ,\\
B&=&\left[r(\om_{-1}+1)+r(-\om_0)\right]\cdot\frac{\mu_{-1}(\om_{-1}+1)\,\mu_0(\om_0-1)}{\mu_{-1}(\om_{-1})\,\mu_0(\om_0)}-\\
&-&r(\om_{-1})-r(-\om_0)-r(\om_0+1)+r(\om_0)\ \ ,\\
C&=&\left[r(-\om_1+1)-r(-\om_1)\right]\cdot\frac{\mu_{-1}(\om_0+1)\,\mu_0(\om_1-1)}{\mu_{-1}(\om_0)\,\mu_0(\om_1)}\cdot\prod_{j\in\Zb}\frac{\mu_{j-1}(\om_j)}{\mu_j(\om_j)}\ \ ,\\
D&=&\left[r(\om_{-1}+1)-r(\om_{-1})\right]\cdot\prod_{j\in\Zb}\frac{\mu_{j+1}(\om_j)}{\mu_j(\om_j)}\ \ .
\end{eqnarray*}
We eliminate the expressions $\mu_k$ for all $k\in\Zb$ with the use of \eqref{eq:ratafelt} and \eqref{eq:miap}
\begin{eqnarray*}
r(z)\cdot\frac{\mu_k(z)}{\mu_k(z-1)}&=&\e{\te_k}\ \ \text{and}\\
r(-z)\cdot\frac{\mu_k(z)}{\mu_k(z+1)}&=&\e{-\te_k}
\end{eqnarray*}
to obtain 
\begin{eqnarray*}
A&=&\sum_{i\ne -1}\Bigl[\e{\te_i-\te_{i+1}}\,r(\om_{i+1})+\e{\te_i-\te_{i+1}}\,r(-\om_i)-r(\om_i)-r(-\om_{i+1})\Bigr]\ \ ,\\
B&=&\e{\te_{-1}-\te_0}\,r(\om_0)+\e{\te_{-1}-\te_0}\,r(-\om_{-1})\,\frac{r(\om_0)}{r(\om_0+1)}-\\
&-&r(\om_{-1})-r(-\om_0)-r(\om_0+1)+r(\om_0)\ \ ,\\
C&=&r(-\om_0)\left(1-\frac{r(\om_1)}{r(\om_1+1)}\right)\e{\te_{-1}-\te_0}\prod_{j\in\Zb}\e{(\te_{j-1}-\te_j)\,\om_j}\frac{Z(\te_j)}{Z(\te_{j-1})}\ \ ,\\
D&=&\left[r(\om_{-1}+1)-r(\om_{-1})\right]\prod_{j\in\Zb}\e{(\te_{j+1}-\te_j)\,\om_j}\frac{Z(\te_j)}{Z(\te_{j+1})}\ \ .
\end{eqnarray*}
For $a<-1,\ b>1$ fixed, let us consider bounded cylinder functions $\vp$, which depend on the variables $\om_a,\,\om_{a+1},\,\dots,\,\om_b$. By stationarity of $\mu^{(\un\te)}$, \eqref{eq:foell} is satisfied for all of them. Hence it is necessary that $A+B+C+D$ does not depend on the variables $\om_a,\,\om_{a+1},\,\dots,\,\om_b$ and its mean is zero according to $\mu^{(\un\te)}$. Only $C,\ D$, and the second term in $B$ are the terms which can contain product of functions of different variables $\om_k$. Each of them is positive by monotonicity of $r$. Thus it follows that each of these three terms must not depend on more than one variable. This implies
\[
\frac{r(z)}{r(z+1)}=\text{constant}=r(0)^2
\]
due to the form of the second term in $B$. The value $r(0)^2$ of this constant is a consequence of \eqref{eq:ratafelt}. Thus we conclude that $r$ is necessarily exponential, the rates are that of the EBL model \eqref{eq:eblr}. $C$ and $D$ can also contain at most one variable, hence we obtain $\te_k=\te_{-1}$ for $k\le-1$ and $\te_k=\te_0$ for $k\ge0$. This means that we have at most two kinds of marginals of $\mu^{(\un\te)}$, one on the left-hand side of the defect tracer and an other one on its right-hand side. In \eqref{eq:foell}, computing the expectation value of $\vp$ times the terms of $A$, summed up for $i\le a-2$ and for $i\ge b+1$, gives zero. The reason for this is that the variables in these terms are independent of the variables which $\vp$ depends on. For the rest of the indices, note that we have a telescopic sum for $A$. Due to this and using the rates of the EBL model, we can simplify our expressions to
\begin{eqnarray*}
A&=&r(-\om_{a-1})-r(\om_{a-1})+r(\om_{b+1})-r(-\om_{b+1})+\\
&+&r(\om_{-1})-r(-\om_{-1})+r(-\om_0)-r(\om_0)\ \ ,\\
B&=&\e{\te_{-1}-\te_0}\,r(\om_0)+\e{\te_{-1}-\te_0-\be}r(-\om_{-1})-\\
&-&r(\om_{-1})-r(-\om_0)-\e{\be}r(\om_0)+r(\om_0)\ \ ,\\
C&=&\left(1-\e{-\be}\right)r(-\om_0)\,\e{(\te_{-1}-\te_0)\,(\om_0+1)}\frac{Z(\te_0)}{Z(\te_{-1})}\ \ ,\\
D&=&\left(\e{\be}-1\right)r(\om_{-1})\,\e{(\te_0-\te_{-1})\,\om_{-1}}\frac{Z(\te_{-1})}{Z(\te_0)}\ \ .
\end{eqnarray*}
It is easy to check that simply choosing $\te_{-1}=\te_0$ does not eliminate the variables $\om_{-1},\,\om_0$ from $A+B+C+D$. Hence this can not be a solution for $\un\te$ to make equation \eqref{eq:foell} be satisfied for all $\vp$. This means that the marginals on the left-hand side of the defect tracer are different from those on the right-hand side. When taking expectation value in \eqref{eq:foell}, this leads to having constant times $\vp$ from the terms containing $\om_{a-1}$ and $\om_{b+1}$ in the first part of the expression of $A$. In order to make \eqref{eq:foell} be satisfied for all $\vp$, it is necessary that we obtain other constants to have zero together with. They can only come from $C$ and $D$. Thus we conclude 
\begin{equation}
\e{\te_{-1}-\te_0}=\e{\be}\label{eq:ered}
\end{equation}
with the use of the form \eqref{eq:eblr} of $r$. In view of \eqref{eq:origmu}, we have the measures 
\begin{equation}
\mu_i(z)=\frac{\e{-\frac{\be}{2}\left(z-\frac{\te_i}{\be}\right)^2}}{\e{-\frac{\te_i^2}{2\be}}Z(\te_i)}=\frac{\e{-\frac{\be}{2}(z-m_i)^2}}{\widetilde{Z}(\be,\,m_i)}\label{eq:mui}
\end{equation}
with $m_i:=\te_i/\be$. We know that the normalization $\widetilde{Z}(\be,\,m)$ in the right-hand side of \eqref{eq:mui} is periodic in the parameter $m$ with period one, which tells us
\[
\frac{Z(\te_{-1})}{Z(\te_0)}=\frac{Z(\te_0+\be)}{Z(\te_0)}=\frac{\e{\frac{(\te_0+\be)^2}{2\be}}\widetilde{Z}(\be,\,\frac{\te_0}{\be}+1)}{\e{\frac{(\te_0)^2}{2\be}}\widetilde{Z}(\be,\,\frac{\te_0}{\be})}=\e{\frac{\be}{2}+\te_0}\ \ .
\]
Using this result together with \eqref{eq:ered} and with the property $\Ev^{(\un\te)}r(\pm\om_i)=\e{\pm\te_i}$, we see that \eqref{eq:foell} is satisfied, which completes the proof.
\end{proof}

The form of the measure described in this theorem shows that the discrete normal distribution of $\om_i,\ i\le-1$ is shifted by $+1$ compared to the distribution of $\om_j,\ j\ge0$. This gives us the picture of a (random) valley with the (randomly) moving defect tracer in its center. Since the position of the defect tracer is not deterministic, we do not see the sharp change between the distribution of the two sides of this valley, if looking the model from outside.

\nemfej{Acknowledgement}
The author is grateful to B\'alint T\'oth for introducing bricklayers' models together with the problem, and for helping him in many questions. Correcting mistakes is also acknowledged to him, as well as to Benedek Valk\'o.
\nocite{ips} 
\nocite{shock}
\nocite{cou}
\nocite{balint}
\nocite{smo}
\nocite{sps}
\nocite{mkps}
\nocite{imes}
\nocite{rez}
\nocite{couse}
\nocite{stochi}

\bibliography{eredeti}
\bibliographystyle{plain}
\end{spacing}
\bigskip\bigskip\bigskip\bigskip
\begin{flushright}
{\renewcommand\arraystretch{0.5}
{\small\begin{tabular}{l}
M\'arton Bal\'azs\\
(balazs@math.bme.hu)\\
Institute of Mathematics,\\
Technical University Budapest\\
1111 Egry J\'ozsef u.\ 1 H \'ep.\ V.\\
Budapest, Hungary
\end{tabular}}}
\end{flushright}
\end{document}